\documentclass[12pt]{amsart}
\usepackage[usenames,dvipsnames]{xcolor}
\usepackage{tikz}
\usepackage{fancyhdr}
\usepackage{txfonts}
\usepackage{graphicx}
\usepackage{epsfig}
\usepackage{mathrsfs}
\usepackage{amssymb}
\usepackage{latexsym}
\usepackage{amsmath}
\usepackage{amsfonts}
\usepackage{amsbsy}
\usepackage{amscd}
\usepackage{subfigure}
\usepackage{verbatim}
\usepackage{color,tikz}
\usetikzlibrary{calc}

\usepackage{amsthm}
\usepackage{amsfonts}
\usepackage{amsbsy,calc}
 \usepackage{graphicx,ifthen,cite}

\newcommand{\D}{{\mathcal D}}

\newcommand{\R}{{\mathbb R}}
\newcommand{\Z}{{\mathbb Z}}
\newcommand{\C}{{\mathbb C}}

\newtheorem{theo}{Theorem}[section]
\newtheorem{defi}{Definition}[section]
\newtheorem{lem}[defi]{Lemma}

\newtheorem{coro}[defi]{Corollary}

\theoremstyle{definition}

\newtheorem{exam}[defi]{Example}
\newtheorem{rem}[defi]{Remark}
\newtheorem{conj}[defi]{Conjecture}
\numberwithin{equation}{section}

\newcommand{\dev}[2][]{{\,\textup{d}}^{#1}#2}

\def\dmu{\dev\mu}

\newif\ifdraft
\drafttrue
\ifdraft
\else
\fi
\numberwithin{equation}{section}
\title{Characterization of self-affine tile digit sets on $\R^n$}
\author{Qian Li,\,Hui Rao}
\date{ 2024}

\begin{document}

\maketitle

\begin{abstract}
  Let $R$ be an $n\times n$ expanding matrix with integral entries. A fundamental question in the fractal tiling theory is to understand the structure of the digit set $\D\subset\Z^n$ so that the integral self-affine set $T(R,\D)$ is a translational tile on $\R^n$. In this paper,
   we introduce a notion of  skew-product-form digit set which is a very general class of tile digit sets. Especially, we show that in the one-dimensional case, if $T(b,\D)$
    is a self-similar tile, then there exists $m\geq 1$ such that
    $$
    \D_m=\D+b\D+\cdots+b^{m-1}\D
    $$
     is a skew-product-form digit set. Notice that $T(b,\D)=T(b^m, \D_m)$, in some sense,
     we completely characterize the self-similar tiles  in $\R^1$. As an application, we establish that all self-similar tiles $T(b,\D)$  where $b=p^{\alpha}q^{\beta}$ contains at most two prime factors are spectral sets in $\R^1$.
\end{abstract}

\bigskip
\section{Introduction}
\bigskip

Let $R$ be an $n\times n$ expanding integer matrix, that is, all the moduli of eigenvalues of $R$ are strictly larger than one; denote $N=|\det(R)|$. Let $\D\subset\Z^n$ be a  set
 of cardinality $\#\D=N$, and we write $\D=\{d_0,\dots, d_{N-1}\}$. Consider the affine iterated function system (IFS) defined by $\{\phi_i\}_{i=0}^{N-1}$, where
$$\phi_i(x)=R^{-1}(x+d_i),\quad 0\le i\le N-1.$$
  The  unique non-empty compact set $T:=T(R,\D)$ satisfying the set-valued functional equation
\begin{equation}\label{eq-32}
  T=\bigcup_{i=0}^{N-1}\phi_i(T)
\end{equation}
  is called the {\it attractor} of the  IFS $\{\phi_i\}_{i=0}^{N-1}$, see \cite{Hu81}.

It is well known (\cite{B1991}) that if $T$ has non-empty interior, then $T$ is a {\it translational tile} , i.e., there exists a discrete set $\mathcal{J}\subset\R^n$ such that
$$T+\mathcal{J}=\R^n\,\,\,\text{and}\,\,\,(T+t_1)^{\circ}\cap(T+t_2)^{\circ}=\emptyset,\,\,t_1\neq t_2\in\mathcal{J}.$$
In that case   $T$ is called an {\it (integral) self-affine tile},  the  set $\D$ is called a {\it  tile digit set} with resect to $R$, and $\mathcal{J}$ is called a {\it tiling set} for $T$. These tiles are referred to as {\it fractal tiles} because their boundaries are usually fractals. There is a large literature on this class of tiles, for instance, \cite{LW1996-1,LW1996-2,LW1996-3,LW1997,GH1994,GH1992,GH1992,SW1999,HLR2003,KL2000,LL2007,GY2006}.

In the theory of self-affine tiles, a fundamental question is to characterize the tile digit sets $\D$   for a given matrix $R$.
This turns out to be a very challenging problem even in $\R^1$. The most basic result is due to Bandt \cite{B1991}, he proved that if $\D$ is a complete residue set modulo $R$, then $T(R,\D)$ is a tile. Such  digit set $\D$ is called  a \textit{standard digit} set in \cite{LW1996-3}. In particular, when $R=p$ is a prime, Kenyon showed that such standard digit sets characterize the tile digit sets \cite{K1991}.
 Another important class of tile digit sets is the \textit{product-form digit sets}, first introduced by Odlyzko \cite{O1978} and later on extended by Lagarias and Wang \cite{LW1996-3}.
   In \cite{LW1996-3}, a modification of the product-form digit set is used to characterize the tile digit sets in $\R$ for $b=p^{\alpha}$, a prime power.
   Lau and Rao  \cite{LR2003} introduce  a \textit{weak product-form digit set} which   was used to characterize the tile digit sets when $b=pq$ is a product of two  primes.  Lai, Lau and Rao \cite{LLR2017} introduced a \textit{modulo product-form  digit set}, and it completely characterized the tile digit sets for $b=p^{\alpha}q$. There are extensions to the higher-dimensional case when $|\det(R)|=p$ is a prime \cite{LW1996-3,HL2004}, and $|\det(R)|=p^2$ a prime square \cite{AL2023}.

In this paper, we introduce  a skew-product-form digit set as follows.

\begin{defi}\label{defi1}
Let $R$ be a $n\times n$ expanding integer matrix, and let $\D\subset\Z^n$ be a   set with $\# \D=|\det R|$.
We call $\mathcal D$ a {\it $m$-stage skew-product-form digit set} with respect to $R$ if there exist $\mathcal{A}=\{a_j\}_{j=0}^{s-1}\subset \Z^n$ and $\mathcal{B}_0,\mathcal{B}_1,\ldots,\mathcal{B}_{s-1}\subset \Z^n$ such that
\begin{equation}\label{eq-41}
\mathcal{D}=\bigcup_{j=0}^{s-1}a_j+R^m \mathcal{B}_j
\end{equation}
and $\mathcal{A}\oplus \mathcal{B}_j$ are complete residue sets modulus $R$, $j=0,1,\ldots, s-1$. In particular, if $m=1$, we call $\D$ a {\it skew-product-form digit set} with respect to $R$.
\end{defi}

\begin{rem}\label{rem-2}
It is worth noting that the structure in \eqref{eq-41} has appeared before.
In the one-dimensional case, An and Lai proved that if $\D$ has the structure in \eqref{eq-41}, then under certain conditions, $T(b,\D)$  is always a spectral set \cite{AL2023}.
In the higher-dimensional case, An and Lau characterized all the tile digit sets when $R=\text{diag}(p,p)$ with $p\geq 2$ being a prime number. In fact, these tile digit sets
 are exactly $k$-stage skew-product-form digit sets \cite{AL2019}.
\end{rem}

 It is a easy matter to show that a skew-product-form digit set is a tile digit set.

\begin{theo}\label{theo1}
If $\D$ is a skew-product-form digit set with respect to $R$, then $T(R,\D)$ is a self-affine tile.
\end{theo}

For $k\ge1$,  denote
\begin{equation}\label{eq-34}
  \D_k=\D+R\D+\cdots+R^{k-1}\D.
\end{equation}

\begin{exam}[\cite{LW1996-3}]\textbf{Product-form.}
If ${\mathcal B_j}=\mathcal{B}$ for $j=0,1,\dots, s-1$ in \eqref{eq-41}, then
$$
\D=\mathcal{A}+R^m{\mathcal B},
$$
and it is a special case of the product-form digit set defined in Lagarias and Wang \cite{LW1996-3}. (A famous example is that
$b=4$ and
$\D=\{0,1,8,9\}.$)  We note that
$$\D_m=(A+RA+\cdots+ R^{m-1}A)+R^m(B+RB+\cdots+R^{m-1}B)$$
is a skew-product-form digit set with respect to $R^m$.

 In general, suppose $\oplus_{j=0}^{m-1} A_j$ is a direct sum decomposition of a complete residue system $\mathcal E$ modulo $R$,
 then
 $$
 \D=\oplus_{j=0}^{m-1} R^jA_j
 $$
 is called a \textit{product-form} according to \cite{LW1996-3}, where $A_j$ can be $\{0\}$.
 Then  $\D_m=U+R^mV$ is a skew-product-form with respect to $R^m$, where
 $$
 U=\oplus_{j=0}^{m-1}\oplus_{i=0}^{m-j-1} R^{i+j}A_j,
 $$
and $V$ contains the other factors in the direct sum decomposition of $\D_m$.
\end{exam}

In general, we show that

\begin{lem}\label{lem3}
If $\D$ is a $m$-stage skew-product-form digit set with respect to $R$, then $\D_m$ is a skew-product-form digit set with respect to $R^m$.
\end{lem}

\begin{exam}[\cite{LR2003}] \textbf{Weak product-form.}
Let
$$
\D \equiv A+R^{m}B \pmod {R^{m+1}},
$$
where $A\oplus B$ is a complete residue set modulus $R$. Then $\D$ is a weak product-form digit set introduced in Lau and Rao \cite{LR2003}. Denote $A=\{a_0,a_1,\ldots,a_{s-1}\}$ and  set
$$B_j=\{u:\,a_j+R^mu\in\D\},\quad 0\le j\le s-1.$$
For any $d\in\D$, there exist $a_j\in A$, $u\in B$ and $x_{j,u}\in\Z^n$ such that
$$d=a_j+R^mu+R^{m+1}x_{j,u}=a_j+R^m(u+Rx_{j,u}).$$
Then
$$\D=\bigcup_{j=0}^{s-1}(a_j+R^mB_j),$$
where $B_j\equiv B\pmod{R}$ for all $j=0,1,\ldots,s-1$. Hence $\D$ is a $m$-stage skew-product-form digit set with respect to $R$.

\end{exam}

\begin{exam}[\cite{LLR2017}] \textbf{Modulo product-form.}
Let $b=12$ and
$$\D=\{0,1,4,8,9,17,25,33,41,72,76,80\}.$$
In \cite{LLR2017},  it is shown that $\D$ is a \textit{modulo product-form}.
Actually, $\D$ is a skew-product-form digit set. Set
$$\mathcal{A}=\{a_j\}_{j=0}^{5}=\{0,1,4,8,9,17\},$$
$$\mathcal{B}_0=\mathcal{B}_2=\mathcal{B}_3=\{0,6\},\,\,\,\mathcal{B}_1=\mathcal{B}_4=\mathcal{B}_5=\{0,2\}.$$
Then
 $\D=\bigcup_{j=0}^{5}(a_j+b\mathcal{B}_j)$
and $\mathcal{A}\oplus\mathcal{B}_j$ are   complete residue systems   module $12$.
\end{exam}

\medskip

Next we consider the converse of Theorem \ref{theo1}. Set $\mathcal{J}_0^*=\Z^n$ and $\mathcal{J}_{k}^*=R\mathcal{J}_{k-1}+\D$ for $k\ge1$.
Clearly  $\{\mathcal{J}_k^*\}_{k=0}^{\infty}$ is a decreasing sequence.
Let $A\subset\R^n$, if there exist a nonsingular integral matrix $B$ and a finite set $\C\subset\Z^n$ such that
$$A=B\Z^n+\C,$$
then we say that $A$ is {\it full-periodic} with period $B$.

\begin{theo}\label{theo2}
Let $R$ be a $n\times n$ expanding integer matrix, and let $0\in\D\subset\Z^n$ be a  set with $\#\D=|\det(R)|$. Suppose that $T(R,\D)$ is a self-affine tile. Then the following are equivalent.

\rm{(i)} There exist an integer $m\ge1$ and a full-periodic set $\mathcal{J}$  with period $R^m$ such that $0\in\mathcal{J}$ and $\mathcal{J}=R^m\mathcal{J}+\D_m$.

\rm{(ii)} There exists an integer $m\ge1$ such that $\D_m$ is a skew-product-form with respect to $R^m$.

\rm{(iii)} There exists an integer $m\ge1$ such that $\mathcal{J}_{m+1}^*=\mathcal{J}_{m}^*$.
\end{theo}

Here we boldly   conjecture that

\begin{conj}\label{rem-1}
If $T(R,\D)$ is a self-affine tile, then there exists an integer $m\ge 1$ such that $\mathcal{J}_{m+1}^*=\mathcal{J}_{m}^*$.
\end{conj}


Instead of characterizing the structure of all digit sets $\D$ with respect to $R$, our focus is on demonstrating that $\D_m$, an iteration of $\D$, exhibits a favorable structure for large values of $m$.
 Usually, to know one `nice' iterated function system of a self-affine tile suffices for  studying the   geometric or analytic properties of  self-affine  tiles.

In particular, in the one-dimensional case,   building upon the result of Lau and Rao \cite{LR2003} regarding the existence, uniqueness and periodicity of the self-replicating tiling, we show that the above conjecture is true and  this essentially  solves the problem of characterizing the tile digit sets $\D$ in $\R^1$  in some sense.

\begin{coro}\label{theo5}
Let $b\ge2$ be an integer and $\D\subset\Z$ be a finite digit set with $0\in\D$ and $\gcd(\D)=1$. Then
$T(b,\D)$ is a self-similar tile if and only if $\D_m$ is a skew-product-form with respect to $b^m$ for some $m\ge1$.
\end{coro}

\medskip

Now we turn to the problem that whether every self-affine tile is a spectral set.
Let $\Omega\subset\R^n$ be a Lebesgue measurable set and let $\mu$ be the Lebesgue measure on $\Omega$ with $0<\mu(\Omega)<\infty$. We call $\mu$ a {\it spectral measure} if there exists a subset $\Lambda$ of $\R$ such that the set of function $\{e^{2\pi i\lambda\cdot x}:\lambda\in\Lambda\}$ forms an {\it orthogonal basis } for $L^2(\mu)$. In this case $\Lambda$ is called a {\it spectrum} for $\mu$ and $\Omega$ is called a {\it spectral set}.

Fuglede initiated the study of general spectral sets in his seminal paper \cite{F1974} and he laid down his famous conjecture stating that $\Omega$ is a spectral set if and only if $\Omega$ is a translation tile. The conjecture is now known to be false in dimension 3 or higher  \cite{KM2006-1, KM2006-2, T2004}. This conjecture remains open in dimension 1 and 2 and the precise classification of spectral sets is still a problem attracting the attention of many researchers. For example, for the case that $\D\subset\Z^{+}$ and $\#\D=4$, Fu, He and Lau \cite{FHL2015} proved that $T(4,\D)$ is a self-similar tile if and only if it is a spectral set.

In \cite{AL2023}, An and Lai provide a sufficient condition for $T(b,\D)$ to be a spectral set. Combining this with Corollary \ref{theo5}, we can settle the ``tiling implies spectrality'' direction of the Fuglede's conjecture in $\R^1$, elying on an auxiliary conjecture proposed by  Coven and Meyerowitz concerning $(T1)$ and $(T2)$ conditions (see Section 3 for the precise definition) \cite{CM1999}.

\medskip

\noindent{\bf C-M conjecture:} Let $A\subset\Z$ be a finite set. If $A$ tiles the integers, then $A$ satisfies $(T1)$ and $(T2)$.

\medskip

\begin{theo}\label{theo16}
Let $b\ge2$ be an integer and $\D\subset\Z$ be a finite digit set with $0\in\D$ and $\gcd(\D)=1$. Suppose the C-M conjecture holds.
If $T(b,\D)$ is a self-similar tile, then $T(b,\D)$ is a spectral set.
\end{theo}

Note that part of this conclusion has already been mentioned in the paper of An and Lai \cite{AL2023},  but we have generalized their findings.
As a corollary of Theorem \ref{theo16}, we have

\begin{coro}\label{theo10}
 Suppose $b=p^{\alpha}q^{\beta}$ contains at most two prime factors. If $T(b,\D)$ is a self-similar tile, then it is a spectral set.
\end{coro}

For the organization of the paper, in section 2, we give the proofs of Theorem \ref{theo1} and Theorem \ref{theo2}. In section 3, we consider the special case in $\R^1$, and we prove Corollary \ref{theo5}, Theorem \ref{theo16} and Corollary \ref{theo10}.

\bigskip

\section{self-affine tile in $\R^n$}

Suppose $T(R,\D)$ is a self-affine tile and $\mathcal{J}$ is a tiling set for it. We say that $\mathcal{J}$ is a {\it self-replicating tiling set}  if
$$\mathcal{J}=R\mathcal{J}+\D.$$

\medskip

First, we  prove  Lemma \ref{lem3}.

\noindent{\bf Proof of Lemma \ref{lem3}.}
Set
$$\D=\bigcup_{i=0}^{s-1}(a_i+R^k\mathcal{B}_i)$$
where $\mathcal{A}=\{a_0,a_1,\ldots,a_{s-1}\}$ and $\mathcal{A}\oplus\mathcal{B}_i$ are complete residue systems  module $R$ for $i=0,1,,\ldots,s-1$. Then
$$\D_k=\bigcup_{i_1=0}^{s-1}\cdots\bigcup_{i_k=0}^{s-1}\left(a_{i_1}+\cdots+
R^{k-1}a_{i_k}+R^k(\mathcal{B}_{i_1}+\cdots+R^{k-1}\mathcal{B}_{i_k})\right).$$
Note that $\mathcal{A}+R\mathcal{A}+\cdots+R^{k-1}\mathcal{A}+\mathcal{B}_{i_1}+\cdots+R^{k-1}\mathcal{B}_{i_k}$ are complete residue sets modulus $R^k$ for all $i_j\in\{0,1,\dots,s-1\}$, $j\in\{1,2,\ldots,k\}$. Then $\D_k$ is a skew-product-form with respect to $R^k$.
\qed

\medskip

Secondly, we prove Theorem \ref{theo1}.

\begin{defi}\label{defi2}\rm{
The {\it upper density} of $\mathcal{J}$ in $\Z^n$ is defined as
$$\overline{D}({\mathcal J})=\varlimsup_{l\to\infty}\frac{1}{(2l)^n}\#\left(\mathcal{J}\cap[-l,l]^n\right).$$
Similarly, we can define the {\it lower density} of $\mathcal{J}$ in $\Z^n$. When the upper density of  $\mathcal{J}$ is equal to the lower one, we refer to the common value as the {\it density} of $\mathcal{J}$, denoted as $D(\mathcal{J})$.
}
\end{defi}

\begin{lem}[\cite{LW1996-2}]\label{lem1}
Let $R$ be an $n\times n$ expanding integer matrix, let $\D\subset\Z^n$ be a   set with $\#\D=|\det(R)|$. Then $T(R,\D)$ is a self-affine tile if and only if $\#\D_k=|\det(R)|^k$ for each $k\ge1$.
\end{lem}

\noindent{\bf Proof of Theorem \ref{theo1}.} Denote $N=\#\D$.
By definition, we can write
$$\mathcal{D}=\bigcup_{j=0}^{s-1}a_j+R\mathcal{B}_j$$
where by setting $\mathcal{A}=\{a_j\}_{j=0}^{s-1}$, $\mathcal{A}\oplus \mathcal{B}_j$ is a complete residue set modulus $R$ for $0\le j\le s-1$.  Let $\mathcal{J}=\mathcal{A}+R\Z^n$. It follows that
\begin{eqnarray*}
  R\mathcal{J}+\D &=& R\mathcal{A}+R^{2}\Z^n+\bigcup_{j=0}^{s-1}(a_j+R\mathcal{B}_j) \\
   &=& \bigcup_{j=0}^{s-1}\left(a_j+R \mathcal{B}_j+R\mathcal{A}+R^{2}\Z^n\right) \\
   &=& \bigcup_{j=0}^{s-1}\left(a_j+R (\mathcal{B}_j+\mathcal{A}+R\Z^n)\right) \\
   &=& \mathcal{A}+R\Z^n.
\end{eqnarray*}
Then $\mathcal{J}=R\mathcal{J}+\D$.

Suppose on the contrary that $N'=\#\D_m\leq N^m-1$ for some $m\geq 2$. Then
$\#\D_{km}\leq (N')^k$. From
$${\mathcal J}=R^{km}{\mathcal J}+\D_{km}=R^{km+1}\Z^n+(R^{km}{\mathcal A}+\D_{km}),$$
we conclude that $D({\mathcal J})=\lim_{k\to \infty} (s\#(\D_{km}))/N^{km+1}=0$, which is a contradiction.
\qed

\bigskip
\bigskip

Thirdly, we  prove Theorem \ref{theo2}.

\noindent{\bf Proof of Theorem \ref{theo2}.} $(\textup{i})\Rightarrow(\textup{ii})$ Suppose that there is a full-periodic set $\mathcal{J}$  with period $R^m$ such that $0\in\mathcal{J}$ and
\begin{equation}\label{eq-5}
  \mathcal{J}=R^m\mathcal{J}+\D_m.
\end{equation}
Denote $N=\#\D$ and $\D_m^{\ast}\equiv\D_m\pmod{R^m}$. Then \eqref{eq-5} implies
\begin{equation}\label{eq-6}
  \mathcal{J}\subset R^m\Z^n+\D_m = R^m\Z^n+\D_m^{\ast}.
\end{equation}
 Since $0\in\mathcal{J}$,  \eqref{eq-5} also implies $\D_m\subset\mathcal{J}$. This together with the fact that $R^m$ is a period of $\mathcal{J}$ implies
\begin{equation}\label{eq-7}
  \D_m^{\ast}+R^m\Z^n\subset\mathcal{J}.
\end{equation}
From \eqref{eq-6} and \eqref{eq-7} we obtain that
\begin{equation}\label{eq-8}
  \mathcal{J}=\D_m^{\ast}+R^m\Z^n.
\end{equation}
Applying \eqref{eq-5} again,
\begin{eqnarray}\label{eq-9}
  \mathcal{J} &=& R^m (\D_m^{\ast}+R^m\Z^n)+\D_m \nonumber\\
   &=& R^m\D_m^{\ast}+\D_m+R^{2m}\Z^n.
\end{eqnarray}

Write $\D_m^{\ast}=\{a_0,a_1,\ldots,a_{s-1}\}$ and set
$$\quad B_j=\{u:\,a_j+R^mu\in\D_m\}, \, \quad 0\le j\le s-1.$$
Then
\begin{equation}\label{eq-11}
  \D_m=\bigcup_{j=0}^{s-1}(a_j+R^m B_j).
\end{equation}
By \eqref{eq-9}, we have
\begin{eqnarray}\label{eq-12}
 {\mathcal J}= \D_m+R^m\D_m^{\ast}+R^{2m}\Z^n &=&\bigcup_{j=0}^{s-1}(a_j+R^m B_j)+R^m\D_m^{\ast}+R^{2m}\Z^n  \nonumber\\
   &=& \bigcup_{j=0}^{s-1}(a_j+R^m(B_j+\D_m^{\ast}+R^m\Z^n));
\end{eqnarray}
by \eqref{eq-8}, we have
\begin{equation}\label{eq-13}
 {\mathcal J}= \D_m^{\ast}+R^m\Z^n=\bigcup_{j=0}^{s-1}(a_j+R^m\Z^n).
\end{equation}
Comparing \eqref{eq-12} and \eqref{eq-13}, we obtain that
 $$B_j+\D_m^{\ast}+R^m\Z^n=\Z^n, \quad j=0,1,\ldots,s-1,$$
 which means that $B_j+\D_m^{\ast}$ contains a complete residue set modulo $R^m$ for each $j$.

Now we claim that $\D_m+R^m\D_m^{\ast}=\D_m\oplus R^m\D_m^{\ast}$ is a direct sum. By \eqref{eq-8}, we have that
\begin{equation}\label{eq-14}
  D(\mathcal J)=\frac{s}{N^m}.
\end{equation}
By \eqref{eq-9}, we have that
\begin{equation}\label{eq-15}
 D(\mathcal J)\leq
\frac{\#(\D_m+R^m\D_m^{\ast}) }{N^{2m}}.
\end{equation}
 To guarantee the equality in the above formula, we must have that
 $\D_m+R^m\D_m^{\ast}=\D_m\oplus R^m\D_m^{\ast}$. The claim  is proved.

Now fix a $j\in \{0,1,\dots, s\}$. From $(a_j+R^m B_j)+R^m\D_m^{\ast}$
is a direct sum we deduce that
  $B_j+\D_m^{\ast}=B_j\oplus\D_m^{\ast}$  is a direct sum. Therefore,
$$
\#B_j\geq N^m/s,
$$
and the equality holds if and only if  $B_j+\D_m^{\ast}$  is a complete residue set modulo $R^m$.
From the decomposition of $\D_m$, we have
$$
N^m=\sum_{j=0}^{m-1} \#B_j\geq s\cdot N^m/s=N^m.
$$
It follows that  every $B_j+\D_m^{\ast}$  is a complete residue set modulo $R^m$ and hence
$\D_m$ is a skew-product-form with respect to $R^m$.

$(\textup{ii})\Rightarrow(\textup{iii})$  Since $\D_m$ is a skew-product-form with respect to $R^m$, we can write
$$\mathcal{D}_m=\bigcup_{i=0}^{s-1}a_i+R^m \mathcal{B}_i$$
where $\mathcal{A}=\{a_i\}_{i=0}^{s-1}$ and $\mathcal{A}\oplus \mathcal{B}_i$ is a complete residue set modulus $R^m$ for $0\le i\le s-1$. Denote $N=\#\D=|\det(R)|$. Then $\mathcal{J}_{m}^{\ast}=R^m\Z^n+\D_m$ and
$$D(\mathcal{J}_{m}^{\ast})=\frac{s}{N^m}.$$
Note that
$$\D_{2m}=\bigcup_{i_1=0}^{s-1}\bigcup_{i_2=0}^{s-1}\left(a_{i_1}+R^ma_{i_2}+R^m (\mathcal{B}_{i_1}+R^m\mathcal{B}_{i_2})\right).$$
It follows that
$$D(\mathcal{J}_{2m}^{\ast})=\frac{sN^m}{N^{2m}}=\frac{s}{N^m}.$$
Since $\{{\mathcal J}_m^\ast\}_{m\geq 0}$ is decreasing, we have
$$D(\mathcal{J}_{m}^{\ast})\le D(\mathcal{J}_{m+1}^{\ast})\le\cdots\le D(\mathcal{J}_{2m}^{\ast})=D(\mathcal{J}_{m}^{\ast}).$$
This forces $D(\mathcal{J}_{m}^{\ast})=D(\mathcal{J}_{m+1}^{\ast})=\cdots=D(\mathcal{J}_{2m}^{\ast})$.
Therefore, since $\mathcal{J}_{m+1}^{\ast}$ is fully periodic, we conclude that
  $\mathcal{J}_{m+1}^{\ast}=\mathcal{J}_m^{\ast}$.

$(\textup{iii})\Rightarrow(\textup{i})$
Set $\mathcal{J}=\mathcal{J}_{m}^{\ast}$. Then $(\textup{iii})$ implies that
$$\mathcal{J}=R\mathcal{J}+\D.$$
It follows that $\mathcal{J}=R^m\mathcal{J}+\D_m$. The theorem is proved.

\qed

\bigskip

\section{The one-dimensional case}
\subsection{Self-similar tile}
\

Now we consider the tiling theory on $\R$. That is, the expanding integer matrix $R=b$ and $\mathcal{D}=\{d_0,d_1,\ldots,d_{b-1}\}$ is a subset of $\R$, where $b\ge2$ is an integer. Then we call $T=T(b,\D)$ a {\it self-similar tile}, and $\mathcal{D}$ a {\it tile digit set} with respect to $b$.


The main result in section 2 still valid in the one-dimensional case.
Furthermore, given the enhanced properties of one-dimensional tiling sets compared to higher-dimensional cases in terms of self-replication and periodicity, we can establish a sufficient and necessary condition for $T(b,\D)$ to be a self-similar tile without making any assumptions, as stated in Theorem \ref{theo2}.

The following lemma verifies the accuracy of Conjecture \ref{rem-1} for the one-dimensional case.

\begin{lem}[\cite{LR2003}]\label{lem-LR2003}
 Assume $0\in\D\subset\Z$ and $\gcd(\D)=1$. If $T(b,\D)$ is a self-similar tile, then there exists a unique self-replicating traslation set $\mathcal{J}$ of $T(b,\D)$ with the property that $\mathcal{J}\subset\Z$. The set $\mathcal{J}$ contains $0$ and is periodic with period $b^m$ for some $m$. More precisely, $\mathcal{J}=\D_m+b^m\Z$.
\end{lem}


\noindent{\bf Proof of Corollary \ref{theo5}.}
The sufficiency immediately follows from Theorem \ref{theo1}. The necessity follows directly from Theorem \ref{theo2} and Lemma \ref{lem-LR2003}.
\qed


%

\bigskip

\subsection{Spectral sets}
\

In this subsection, we aim at proving ``tiling implies spectrality''.
Before this, we first recall the results of  Coven  and Meyerowitz \cite{CM1999} concerning
 direct sum decompositions of  cyclic groups. Let $A$ be a finite set of integers. Then $P_A(x)=\sum_{a\in A}x^a$ is a polynomial such that $\# A=P_A(1)$. Let $\Phi_s(x)$ be the $s$th cyclotomic polynomial, which is the minimal polynomial of $e^{\frac{2\pi i}{s}}$ in $\Z[x]$. Let $S_A$ be the set of prime powers $s$ such that the $s$th cyclotomic polynomial $\Phi_s(x)$ divides $P_A(x)$.

Coven and Meyerowitz introduced the following two conditions for $A(x)$ which guarantees tiling for $A$ for some $\Z_N$.

\rm{(T1)} \quad $P_A(1)=\prod_{s\in S_A}\Phi_s(1)$.

\rm{(T2)} \quad For all distinct prime powers $s_1,\ldots,s_k\in S_A$, we have $\Phi_{s_1\cdots s_k}(x)\mid P_A(x)$.

Coven and Meyerowitz proved the following theorem \cite{CM1999}.

\begin{theo}[\cite{CM1999}]\label{theo13}
If $\# A$ has at most two prime factors, then $A$ tiles the integers if and only if $A$ satisfies $(T1)$ and $(T2)$.
\end{theo}

The $(T1)$ and $(T2)$ conditions are also sufficient to obtain a spectrum on $\Z_n$. This is due to {\L}aba.

\begin{theo}[\cite{L2002}]\label{theo14}
Let $A$ be finite set of non-negative integers which satisfying $(T1)$ and $(T2)$, then $A$ is a spectral set in $\Z_N$ with a spectrum
$$\Lambda=\left\{\sum_{s\in S_{A}}\frac{k_s}{s}:\,\,k_s\in\{0,1,\ldots,p-1\}\,\,\textup{if}\,\,s=p^{\alpha}\right\}$$
and $N=\operatorname{lcm}(S_A)$.
\end{theo}

\begin{defi}
We say that $(b,\D,\mathcal{C})$ forms an (integral) Hadamard triple if $\D$, $\mathcal{C}\subset\Z$ and $\delta_{b^{-1}\D}$ is a spectral measure with spectrum $\mathcal{C}$. Equivalently,
$$\frac{1}{\sqrt{\#\D}}\left(e^{2\pi i\frac{dc}{b}}\right)_{d\in\D,c\in\mathcal{C}}$$
is a unitary matrix.
\end{defi}

If $A$ is a finite set of real number, we will write
$$\delta_{A} = \frac1{\#A} \sum_{a\in A}\delta_a$$
where $\delta_a$ is the Dirac measure at $a$.
The Fourier transform of a Borel probability measure is defined to be
$$\hat{\mu}(\xi)=\int e^{-2\pi i\xi x}\dmu(x).$$
Then
$$\hat{\delta}_{A}(\xi)=\frac{1}{\#A}\sum_{a\in A }e^{-2\pi i a\xi}.$$
We recall the following equivalent conditions \cite{LW2002}.

\begin{lem}[\cite{LW2002}]\label{lem2}
The following are equivalent.

\rm{(i)}  $(N,A,L)$ forms a Hadamard triple.

\rm{(ii)} $\delta_{\frac{A}{N}}$ is a spectral measure with spectrum $L$.
\end{lem}

%

In \cite{AL2023}, An and Lai give a sufficient condition for $T(b,\D)$ to be a spectral set.

\begin{theo}[\cite{AL2023}]\label{theo15}
Let $N\ge2$ be an integer and $\mathcal{A}=\{a_j:\,j=0,1,\ldots,s-1\}$ be a subset of integers and for each $j$, we have $\mathcal{B}_j$ as another subset of integers. If there exists $r\ge0$ such that
$$\D=\bigcup_{j=0}^{s-1}(a_j+N^r\mathcal{B}_j)$$
and we have the following conditions for $\mathcal{A}$, $\mathcal{B}_j$, $L_1$ and $L_2$

\rm{(i)} \quad $(N,\mathcal{A},L_1)$ and $(N,\mathcal{B}_j,L_2)$ are Hadamard triples for all $j=0,1,\ldots,s-1$;

\rm{(ii)} \quad $(N,\mathcal{A}\oplus\mathcal{B}_j,L_1\oplus L_2)$ are Hadamard triples for all $j=0,1,\ldots,s-1$.\\
Then $T(N,\D)$ is a spectral set.
\end{theo}

\medskip

\noindent{\bf Proof of Theorem \ref{theo16}.}
Applying Corollary \ref{theo5}, there exists $m\ge1$ such that $\D_m$ is a pseudo-product-form with respect to $b^m$.
Denote
$$\mathcal{D}_m=\bigcup_{j=0}^{s-1}(a_j+b^m \mathcal{B}_j)$$
where $\mathcal{A}=\{a_j\}_{j=0}^{s-1}$ and $\mathcal{A}\oplus \mathcal{B}_j$ is a complete set of residues modulo $b^m$ for $0\le j\le s-1$. Set
$$\Z_{b^m}=\{0,1,\ldots,b^m-1\}.$$
Then
\begin{equation}\label{eq-35}
  \mathcal{A}\oplus \mathcal{B}_j\equiv \Z_{b^m}\pmod{b^m},\quad\forall 0\le j\le s-1.
\end{equation}
 Notice that \eqref{eq-35} yields that $S_{\mathcal{B}_j}$ are all the same for $0\le j\le s-1$.
For any $0\le j\le s-1$, we denote
$$L_1=\sum_{i=1}^{r_1}\frac{1}{p_{i}^{\alpha_i}}\{0,1,\ldots,p_i-1\},\,\,\,
L_2=\sum_{k=1}^{r_2}\frac{1}{q_k^{\beta_k}}\{0,1,\ldots,q_k-1\}$$
where all $p_i,q_i$ are prime numbers and $\alpha_i,\beta_i$ are positive integers.
Applying Theorem \ref{theo14}, we know that $L_1$ and $L_2$ are spectra of $\delta_{\mathcal{A}}$ and $\delta_{\mathcal{B}_j}$, respectively. Therefore, $(N_1,\mathcal{A},N_1L_1)$ and $(N_2,\mathcal{B}_j,N_2L_2)$ are Hadamard triples for $0\le j\le s-1$ by Lemma \ref{lem2}, where $N_1=\operatorname{lcm}(S_A)$ and $N_2=\operatorname{lcm}(S_{B_j})$. Combining with the fact that $b^m=\operatorname{lcm}(S_A,S_{B_j})$, we have $(b^m,\mathcal{A},b^mL_1)$ and $(b^m,\mathcal{B}_j,b^mL_2)$ are Hadamard triples for $0\le j\le s-1$. Note that
\begin{equation}\label{eq-40}
  \# (L_1\oplus L_2)=\prod_{i=1}^{r_1}p_i\prod_{k=1}^{r_2}q_k
  =(\#\mathcal{A})(\#\mathcal{B}_j)=b^m=\operatorname{lcm}(S_{\Z_{b^m}}).
\end{equation}
It follows that the power of each prime power in $S_{\Z_{b^m}}$ appears continuously. Therefore, $b^m(L_1\oplus L_2)$ is a complete set of residues modulo $b^m$, which means $(b^m,\mathcal{A}\oplus\mathcal{B}_j,b^mL_1\oplus b^mL_2)$ are Hadamard triples for all $j=0,1,\ldots,s-1$. Then $T(b^m,\D_m)=T(b,\D)$ is a spectral set by  Theorem \ref{theo15}.
\qed

\medskip

\noindent{\bf Proof of Theorem \ref{theo10}.} The conclusion follows directly from Corollary \ref{theo5}, Theorem \ref{theo13} and Theorem \ref{theo16}.
\qed

\bigskip

\bigskip

\end{document}